\RequirePackage{fix-cm}
\documentclass[twocolumn]{svjour3}          
\smartqed  
\usepackage{graphicx}
\usepackage[misc,geometry]{ifsym}
\usepackage{algorithm}
\usepackage{placeins}
\usepackage{algpseudocode}
\usepackage{subcaption}
\usepackage{amsfonts}
\usepackage{amssymb,amsmath}
\usepackage{todonotes}
 \journalname{Engineering with Computers}
\begin{document}

\title{Mesh Improvement Methodology for 3D Volumes with non-Planar Surfaces
}


\author{Alan Kelly         \and
        Lukasz Kaczmarczyk \and 
	Chris J. Pearce
}

\institute{A. Kelly {\Letter} \and  L. Kaczmarczyk \and  C.J. Pearce  \at
              School of Engineering \\
	      University of Glasgow,\\
	      Glasgow G12 8LT, UK\\
              \email{a.kelly.2@research.gla.ac.uk} \\          
          }

\date{Received: date / Accepted: date}

\maketitle

\begin{abstract}
A mesh improvement methodology is presented which aims to improve the quality of the worst elements in 3D meshes with non-planar surfaces which cannot be improved using traditional methods. A numerical optimisation algorithm, which specifically targets the worst elements in the mesh, but is a smooth function of nodal positions is introduced. A method of moving nodes on curved surfaces whilst maintaining the domain geometry and preserving mesh volume is proposed. This is shown to be very effective at improving meshes for which traditional mesh improvers do not perform well.
\keywords{Mesh Optimisation \and Volume Preservation \and Log-Barrier \and Patch-Improvement}
\end{abstract}

\section{Introduction}
\label{intro}
In the context of the Finite Element Method (FEM), high quality meshes can be crucial to obtaining accurate results. The quality of an element can be described as a numerical measure which estimates the effect that the shape of an element will have on the accuracy of an analysis,~\cite{Shewchuk2002}. It can be shown that poor quality elements can result in both discretisation errors and poor conditioning of the stiffness matrix. In the extreme, a single poor element can render a problem intractable. Therefore, a high quality mesh is crucial to performing an accurate analysis. 

The field of mesh optimisation is complex and has now become an area of research in its own right. Numerical optimisation is the process of maximising or minimising an objective function, subject to constraints on the solution. When this is applied to a finite element mesh it is referred to as mesh optimisation, where the mesh quality is the objective function and the constraints include, for example, the domain geometry and maximum element size. In order to make the process of mesh optimisation more straightforward for the analyst, we aim to create a set of tools which makes it possible to improve complex meshes used in actual simulations, in as simple a manner as possible. In doing this, we are attempting to simplify a very complex process; this paper will explain the problems encountered and the solutions to these problems.
\subsection{Motivation}
The motivation for this project is the need for high quality meshes for problems with evolving geometries, such as fracturing solids, moving fluids and biological materials. For complex three-dimensional geometries, automatic mesh generators do not always create meshes of sufficient quality to ensure a sufficient level of accuracy in the solution. This is further complicated by the need to have an adapting mesh that can resolve the evolving geometry. Issues with the solution of these problems can be traced back to poor quality meshes and although there are a number of tools already available for improving mesh quality, none of them have matched the needs of the authors.

For the creation and evolution of a mesh, the positioning of the nodes is determined by physics of the problem being analysed and therefore any alteration of their positioning must be compatible with the physics of the problem, i.e. the geometry and the volume of the domain must be preserved. Therefore a method of improving mesh quality by moving surface nodes but without changing the geometry or volume of the domain is necessary. Such a method has been developed and is described in Section \ref{surface}.
\subsection{Implementation}
\label{sec1.2}
This work was implemented using the Mesh Quality Improvement Toolkit (Mesquite) as a platform,~\cite{Knupp2012}. The architecture of this library makes it ideal to use as a base for the development and testing of new algorithms. In addition, BLAS and Lapack were used for numerical operations. Both Mesquite's native algorithms~\cite{Knupp2012} and Stellar~\cite{Klingner2008}, a mesh optimisation program, were used to assess the results obtained. Stellar, whilst a very powerful program, has restrictions that mean it has limited application to the kind of problems that motivated this work. The main issues are that the user has very limited control over the optimisation process and that it was developed with the goal of achieving the highest mesh quality possible, regardless of time taken. The user has no control over termination criteria, limited control over what improvement operations are performed and it is not easily integrated into other projects. Many mesh generation packages (e.g. Cubit \cite{CUBIT}) allow for data, boundary conditions for example, to be added to the mesh as part of the generation process. This data will be lost when a mesh is added to Stellar due to the changes in mesh topology and the addition/elimination of nodes, (operations the user has very limited control over). Also, the user cannot fix arbitrary nodes which is often essential in FEA. Although Stellar can improve meshes effectively, the loss of control over the optimisation process renders it unsuitable for use in many FE simulations. However, as an academic package, it is very powerful; it demonstrates the quality which can be achieved through mesh improvement operations and also demonstrates the effectiveness of certain operations. For these reasons Stellar is an ideal tool for comparison of results and for deciding which operations are worth implementing.

\subsection{Optimisation-Based Mesh Smoothing}
\label{sec:1.3}
Mesh smoothing is the process of improving mesh quality without changing the mesh topology~\cite{Freitag}. Mesh topology refers to the nodes of the mesh and the elements which these nodes lie on. There are many existing mesh smoothing algorithms, the most famous of which is Laplacian smoothing. Laplacian smoothing involves moving a vertex to the average of its connected neighbours and is applied to each mesh vertex in sequence and is repeated several times. It has been shown to be somewhat effective with 2D triangular meshes, but is much less effective in 3D~\cite{Klingner2008}. Although Laplacian smoothing is computationally cheap, there is no guarantee of mesh improvement. It is even possible that inverted elements will be created~\cite{Chen2003} when the domain is not convex~\cite{Donea2004}. Much more sophisticated mesh smoothing algorithms have been developed which are based on numerical optimisation techniques. Techniques such as these are referred to as optimisation-based smoothers. These methods require a means of expressing the quality of an element numerically and of combining the qualities of every element in the mesh into a single numerical measure. A numerical measure which effectively captures mesh quality requirements is described in Section \ref{QM}. 

Mesh quality optimisation requires an objective function which combines the qualities of a group of elements into a scalar value. For example, one could express the quality of a mesh as the sum of the qualities of every element. This objective function would then be minimised, or maximised, depending on the choice of quality measure, to improve the quality of the mesh. However, as previously stated, one poor element may render a problem unsolvable. A simple objective function such as the one described above would be very good at improving average element quality but would not improve the worst element, since one poor quality element would not stand out. Such an objective function may even invert some elements, as one negative number may not sufficiently influence the objective function. Therefore, it is desirable to use an objective function that targets the quality of the worst element.

At first glance, an Infinity Norm seems like an ideal objective function. This is where the quality of a group of elements is expressed as the quality of the worst element. In this case, any attempt to optimise the mesh will improve the worst element. However, nodes are shared between elements. So if a node is moved to increase the quality of one element, the quality of adjoining elements may be adversely affected. As the infinity norm contains no information about the adjoining elements' quality, there is no way of knowing when the element being improved is no longer the worst element in the mesh. Therefore, such an objective function is described as being non-smooth. A non-smooth optimisation algorithm was developed by~\cite{Freitag1995}, which enabled the improvement of the worst element in a mesh. This algorithm achieved very high quality results and is utilised in Stellar~\cite{Klingner2008}. This approach works by calculating the search directions for the nodes of the worst element and attempting to predict the distance each node may be moved in this direction until the element is no longer the worst. As element quality is a function of nodal positions, a first order Taylor Series expansion of the quality of every affected element may be used to approximate the point at which the element being improved is no longer the worst.

A genuinely smooth objective function which penalises the worst element in a mesh to such an extent that the improvement process focuses on this element should, in theory, yield better results in a shorter analysis time since the objective function contains information about the quality of all elements, so there is no requirement to approximate the point at which the quality of the worst element changes. An objective function which meets this smoothness criterion and which also adequately penalises the worst element is described in Section \ref{subsubsec:2.1.2}.
\subsection{Meshes}
\begin{figure*}
  \centering
 \hspace{-20mm}\begin{subfigure}[b]{0.35\textwidth}
    \center{\includegraphics[height=6cm]{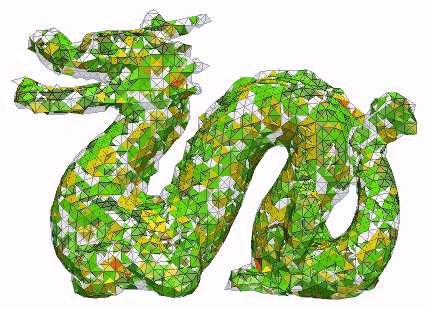}}
    \raggedleft\caption{\textsc{Dragon}, 32959 Tetrahedra \cite{Klingner:2007:ATM}\hspace{5mm}}
    \label{dragon_orig}
  \end{subfigure}
  \hspace{35mm}
 \begin{subfigure}[b]{0.3\textwidth}
    \center{\includegraphics[height=6cm]{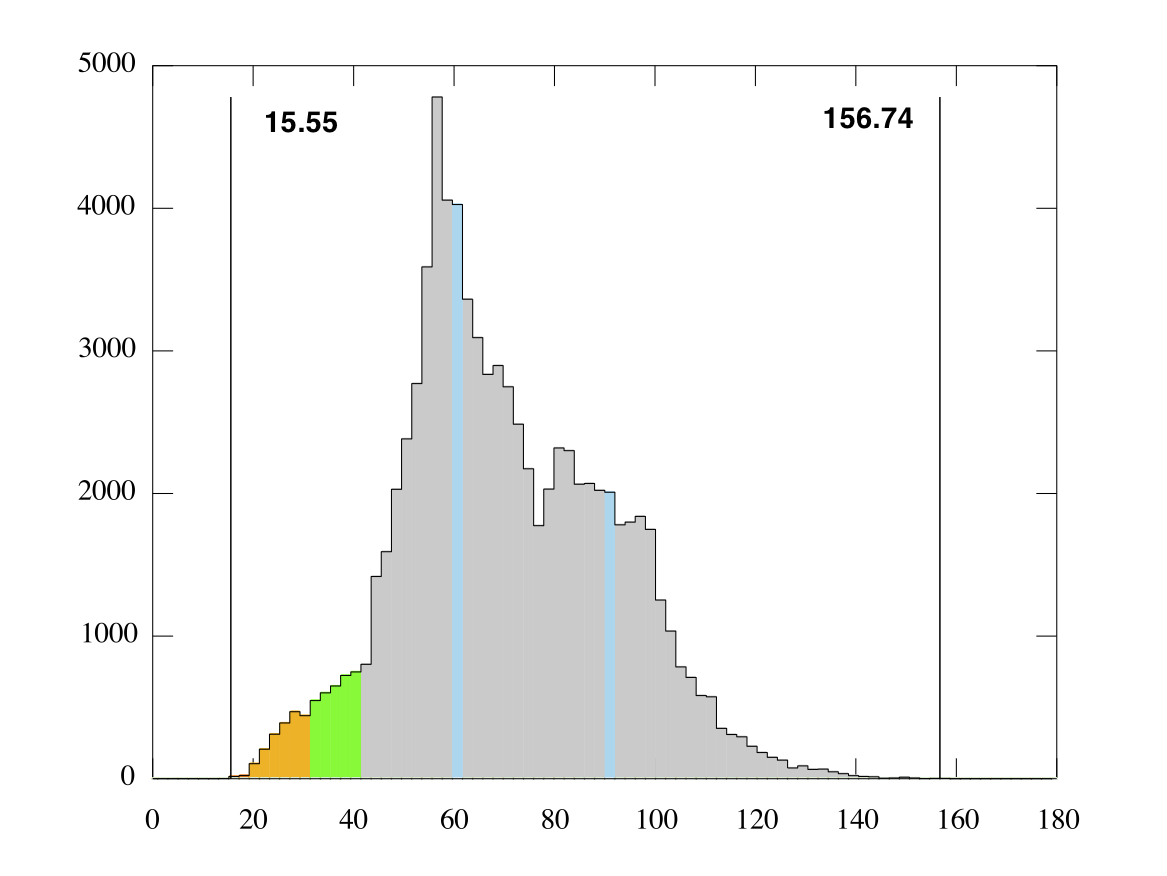}}
    \label{dragon_orig_hist}
  \end{subfigure}\\
    \hspace{-20mm}\begin{subfigure}[b]{0.35\textwidth}
    \center{\includegraphics[height=5cm]{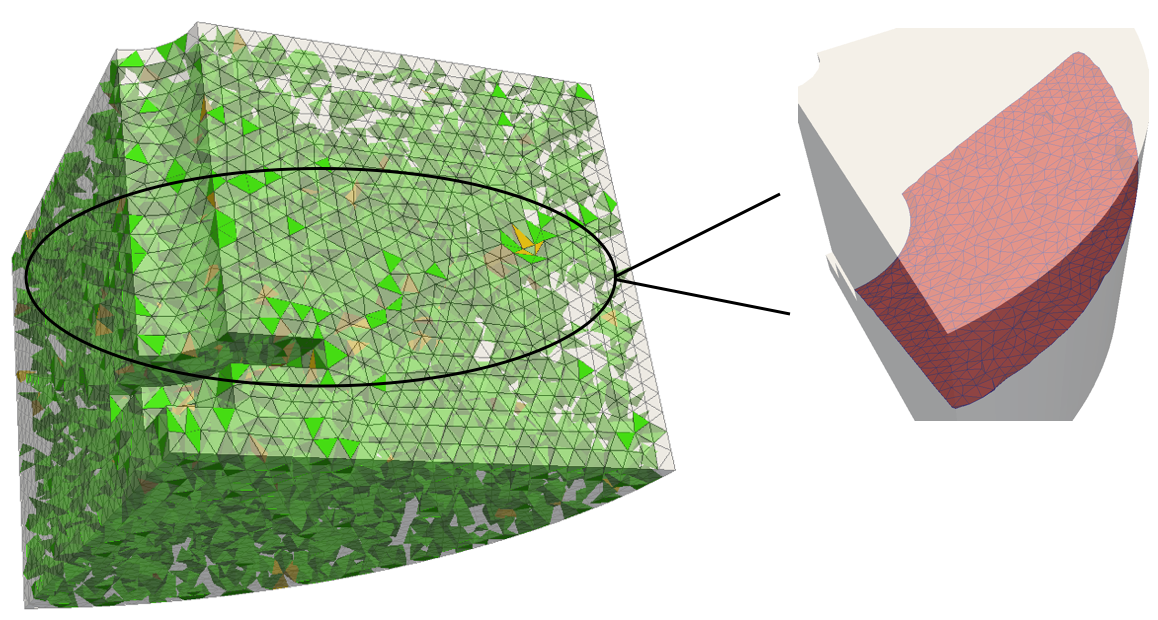}}
    \raggedright\caption{\textsc{Concrete Cylinder}, 73684 Tetrahedra}
    \label{pull_orig}
  \end{subfigure}
   \hspace{35mm}
  \begin{subfigure}[b]{0.3\textwidth}
   \center{\includegraphics[height=6cm]{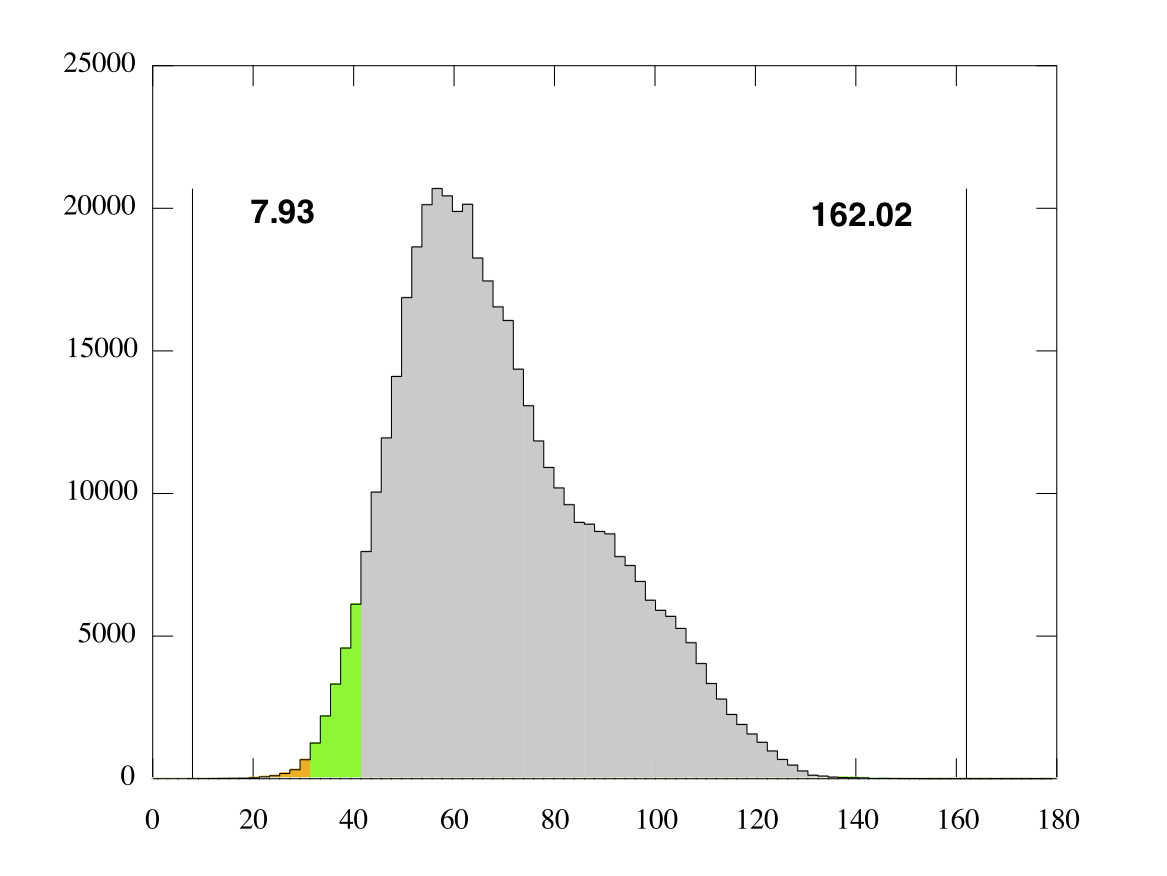}}
    \label{pull_orig_hist}
   \end{subfigure}\\
\hspace{-20mm}\begin{subfigure}[b]{0.35\textwidth}
    \center{\includegraphics[height=6cm]{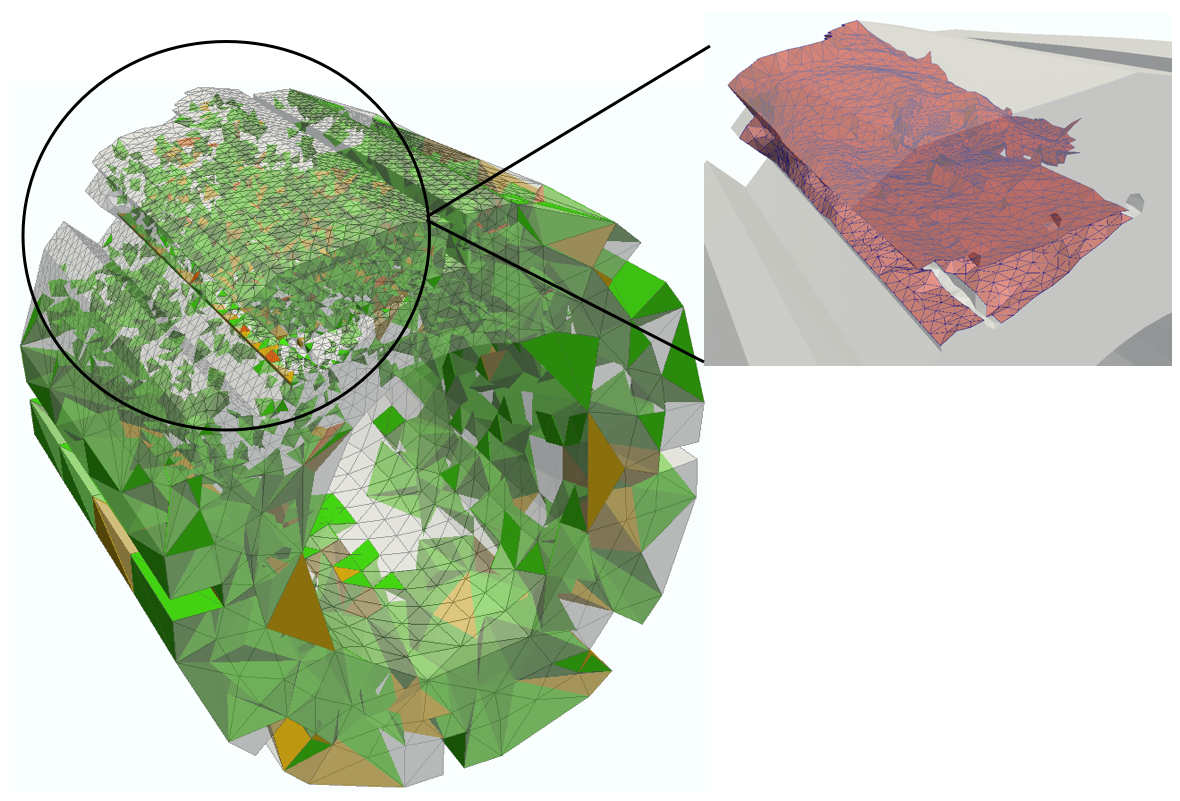}}
    \raggedright\caption{\textsc{Graphite Brick}, 100556 Tetrahedra}
    \label{brick_orig}
  \end{subfigure}
  \hspace{35mm}  
  \begin{subfigure}[b]{0.3\textwidth}
    \center{\includegraphics[height=6cm]{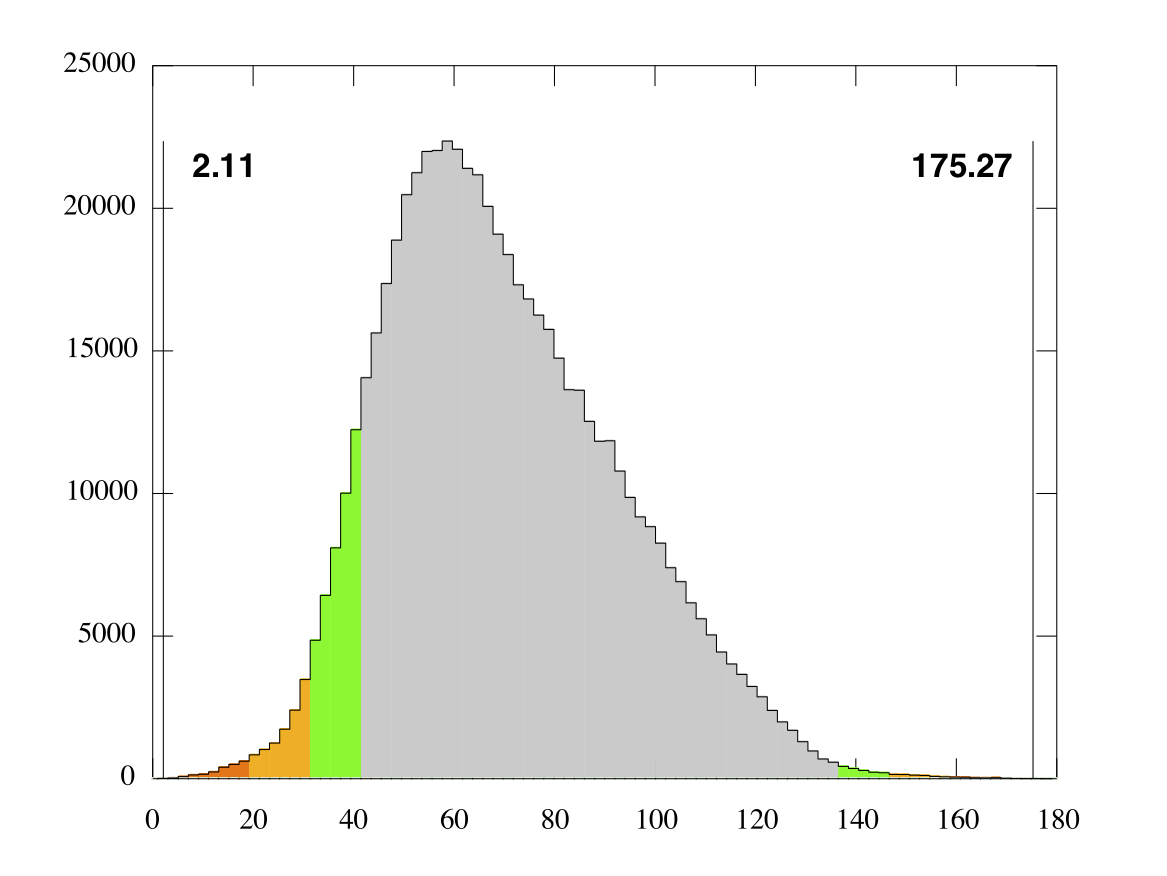}}
    \label{brick_orig_hist}
  \end{subfigure}\\

  \caption{Meshes used for Testing with Histograms of the range of dihedral angles $($the height of blue columns have been divided by 20 due to the many occurrences of these angles$)$. Internal crack surfaces are shown in red. Orange tetrahedra have angles under 20$^\circ$ or greater than 160$^\circ$, yellow tetrahedra have angles between 20$^\circ$ and 30$^\circ$ or 150$^\circ$ and 160$^\circ$ and green tetrahedra have angles between 30$^\circ$ and 40$^\circ$ or 140$^\circ$ and 150$^\circ$.}
  \label{plot:mesh}
  %
  %
\end{figure*}

Three different meshes are used throughout this paper to demonstrate the effectiveness of the algorithms described in the following sections (see Figure \ref{plot:mesh}). \textsc{Dragon} was generated by Isosurface Stuffing, \cite{Labelle:2007:ISF}. Both \textsc{Concrete Cylinder} and \textsc{Graphite Brick} were generated using Cubit \cite{CUBIT}. These meshes were chosen since their complex geometries demonstrate clearly the need for sophisticated mesh optimisation algorithms and the effectiveness of the techniques described in this paper. The crack surface in \textsc{Concrete Cylinder} was formed by simulating a piece of steel which was encased in a concrete cylinder being pulled until the specimen failed. The crack surface in \textsc{Graphite Brick} was formed by simulating an external force being placed on a sample of graphite used in the construction of nuclear power plants. The simulation process is explained in detail in \cite{Luk:2013}. As both crack surfaces are generated by simulating physical phenomena, we don't have any additional information about their smoothness.

\section{Mesh Improvement Methodology}
\label{sec:2}
The focus of this research is to produce practical tools which can be used in a variety of problems. These tools must be easy to use, powerful, efficient and easy to integrate into existing codes. In the case of domains which are constantly evolving, it is often necessary to perform mesh improvement as a continuous process during an analysis. Therefore, algorithms which can quickly improve the worst elements of large meshes are required.
\subsection{Quality Measure}
\label{QM}
Finding a suitable quality measure that provides an accurate estimate of an element's effects in terms of discretisation/interpolation error and stiffness matrix condition is challenging and is a very active area of research in itself. There are many measures in existence and these may be further studied in \cite{Shewchuk2002}. A tetrahedron has many properties which determine its effect on the accuracy of a FEM simulation, however the dihedral angles formed between the faces of a tetrahedron have been shown to be of greatest importance, \cite{Klingner:2007:ATM}. These angles range between 0$^\circ$ and 180$^\circ$. Large dihedral angles have been shown to result in interpolation errors and small dihedral angles affect the conditioning of the stiffness matrix. At first glance, one would think that a quality measure based on the dihedral angles of a tetrahedron would be the most suitable and such measures do exist, for example the minimum sine measure which expresses the quality of an entire mesh as the minimum of the sines of all the dihedral angles in that mesh. Thus, this measure is \textit{non-smooth} meaning that it cannot be used with most mesh optimisation algorithms. Algorithms which can cope with a non-smooth measure are very computationally expensive, thus making them incompatible with the overall project goal of creating practical tools  which may be easily and inexpensively integrated into existing projects. When the smoothing process is a part of a monolithic FE simulation, that is to say smoothing is a part of the overall algorithm, a smooth objective function allows the use of a Newton method which converges quadratically.

To achieve this, it was decided to use a measure called the Volume-Length quality measure. Although this measure does not directly measure poor dihedral angles, it has been shown to be very effective at eliminating such undesirable angles, thus improving stiffness matrix conditioning and interpolation errors \cite{Shewchuk2002},\cite{Klingner2008}. As the Volume-Length measure is a smooth function of vertex positions and its gradient/Hessian are straightforward and computationally cheap to calculate, it is ideal for our requirements. This measure is normalised so that an equilateral element has quality 1 and a degenerate element (zero volume) has quality 0. This measure does not directly correlate with a measure of the dihedral angles; for example, an element with a very low Volume-Length ratio could in fact have good dihedral angles. However, in practice, it has been found that using this measure is very effective. Many other quality measures were considered such as Mesquite's Jacobian based measures such as the Ideal Weight Mean Ratio measure. This measure, in its standard form is very similar to the Volume-Length ratio in that a degenerate element has quality 0, an equilateral has quality 1 and an inverted element has negative quality. Tests indicate a closer correlation between between the Volume-Length ratio and an element's dihedral angles then the Ideal Weight Mean Ratio.
\begin {equation}
  \label{eqn:VL}
  q=6\sqrt2\frac{V}{l^{3}_{rms}}
\end{equation}
\noindent
where $V$ is the volume of a tetrahedral element, $l_\textrm{rms}$ is the root mean square of the element's edge lengths and $q$ is the element quality.

\subsection{Worst Element Improvement Algorithms}
\label{subsec:2.1}

\subsubsection{The Log-Barrier Objective Function}
\label{subsubsec:2.1.2}

This section describes an objective function which both satisfies the smoothness criterion described in Section \ref{sec:1.3} and punishes the worst element in the mesh. This is achieved by expressing the quality of every element as a function of the worst element, equation (\ref{eqn:1}), and is referred to as a log-barrier function. Expressing the quality of a group of elements in this manner ensures that the optimisation process is always focused on the worst element. The Log-Barrier function, its gradient and Hessian are defined as:

\begin{equation}
  \label{eqn:1}
  I=\frac{q^{2}}{2(1-\gamma)}-log(q-\gamma)
\end{equation}

\noindent where $\gamma=b*q_{min}$

\begin{equation}
  \label{eqn:2}
  \mathbf{\mathbf{f}} = (\frac{q}{1-\gamma}-\frac{1}{q-\gamma}){\nabla q}
\end{equation}

\begin{equation}
  \label{eqn:3}
  \mathbf{\mathbf{S}}= {\nabla q}[\frac{1}{1-\gamma}-\frac{1}{(q-\gamma)^2}]{\nabla q^{T}}+[\frac{q}{1-\gamma}-\frac{1}{q-\gamma}]{\nabla^2 q}
\end {equation}
  
\noindent where $q$ is the element quality and $\gamma$ is the barrier which is a function of \emph{b} and the worst element in the mesh, q$_{min}$, ${\nabla q}$ the gradient and ${\nabla^2 q}$ the Hessian of the quality measure.

It has been found that choosing a barrier constant term \emph{b} in the range 0.75-0.95 is most effective. The optimisation process starts with a lower value of \emph{b} and becomes more aggressive with \emph{b} increasing. Smaller values of \emph{b} tend to increase average element quality as the worst elements are not punished as harshly, whereas higher \emph{b} tends to improve the quality of the worst element.

\begin{figure*}
  \centering
  \begin{subfigure}[b]{0.3\textwidth}
    \includegraphics[height=3.5cm]{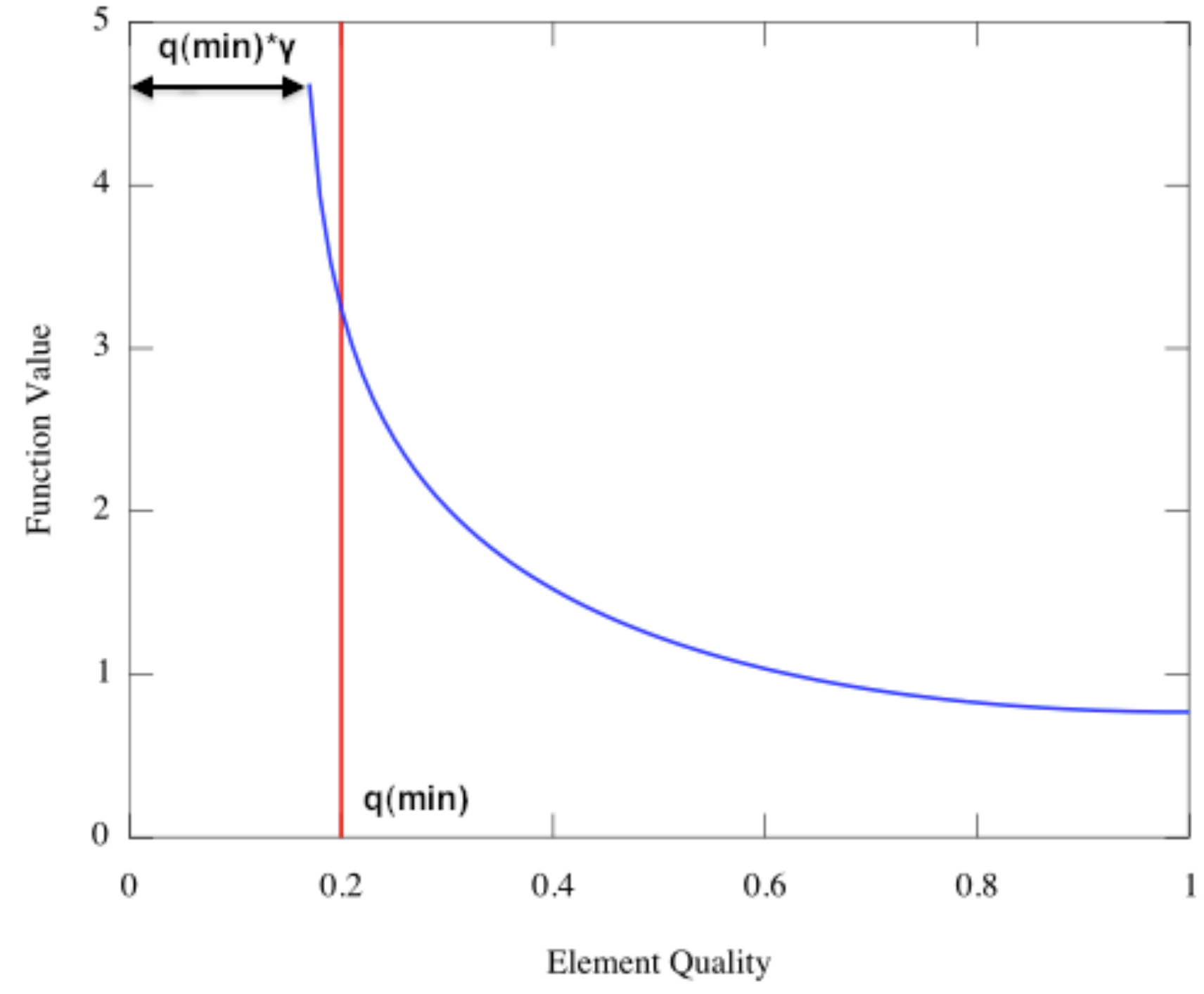}
    \caption{(a) Log-Barrier}
  \end{subfigure}
  \qquad \hspace{25mm}
  \begin{subfigure}[b]{0.3\textwidth}
    \includegraphics[height=3.5cm]{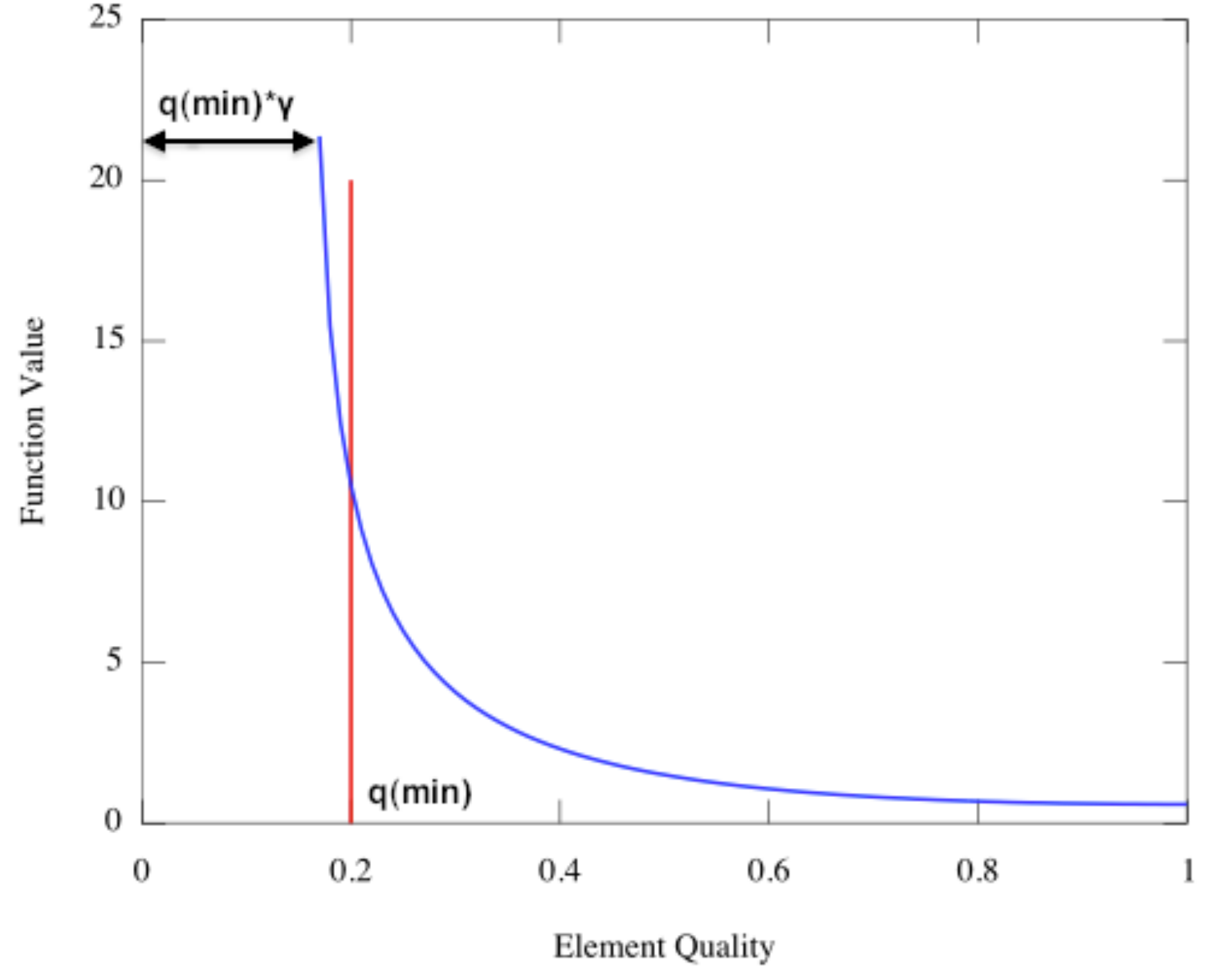}
    \caption{ (b) Log-Barrier Squared}
    \label{barriersquare}
  \end{subfigure}
  \caption{Plot of the Log-Barrier function (equilateral element has quality 1 and a degenerate element (zero volume) has quality 0) (\emph{$\gamma$}=0.8)}
  \label{plot:1}
  %
  %
\end{figure*}

Figure \ref{plot:1} shows the Log-Barrier function graphically. It can be seen that the function rapidly increases as the quality of the element reduces, thus achieving our aim of punishing the worst element. Figure \ref{barriersquare} demonstrates that this effect is further magnified by squaring the Log-Barrier function.

A Newton based solver is used to optimise the mesh, which solves the following nonlinear system of equations to determine nodal positions:

\begin{equation}
  \label{Sdx}
  \mathbf{S}\delta\mathbf{X}=-\mathbf{f}
\end{equation}

This process is repeated several times for the Log-Barrier function as the parameters change as the worst element changes. Unlike traditional mesh optimisation, which is allowed to run until it is deemed to have converged or some other termination criteria has been achieved, Log-Barrier optimisation performs one pass over each patch and then the worst element quality is re-calculated and $\gamma$ is updated to reflect this. This ensures that the optimisation process is always aggressively tackling the worst element.

The Log-Barrier function also has several other very useful features. It comes with an invertibility guarantee - if the initial mesh is valid, that is to say a mesh without any inverted or negative volume elements, is input, no inverted elements will be created. If the initial mesh is invalid, the Log-Barrier function can untangle it as the quality is always chosen to be worse than the worst element. The form of the Log-Barrier function adopted here is different to the form described by \cite{Sastry2011}. The differences between both methods should be investigated to compare their respective merits.

\subsection{Comparison of Measures}
\label{compmeas}
Mesh optimisation was performed on the three meshes shown in Figure \ref{plot:mesh}.  The Log-Barrier function combined with the Volume-Length quality measure was compared with both Mesquite's Ideal Weight Inverse Mean Ratio quality measure combined with an Infinity Norm objective function (found to be the most effective objective function in Mesquite at improving the worst element in a mesh) and Stellar. Each smoothing algorithm was run until convergence, with no restrictions on time, to measure the highest quality each smoother could achieve. The results may be seen in Table \ref{tab:1}. For all three meshes, the Log-Barrier function achieved greater improvement than either Mesquite's native algorithms or Stellar's smoothing algorithms. It must be noted that, although all topological transformation functionality was turned off in Stellar, some changes were made to the mesh connectivity in all three cases, including the deletion of nodes. The authors believe that the combination of the Log-Barrier function with a suitable quality measure, such as the Volume-Length measure, is the most effective method of optimising mesh. Three notes must be made about these results at this point. First, for all three smoothers, surface vertices were unrestrained, meaning that geometry and volume were not preserved. Second, only mesh smoothing algorithms are compared here. Third, Stellar can achieve significantly better results  when all of its functionality is enabled. However, this often results in a mesh with significantly fewer vertices than the original mesh. This last point demonstrates the potential improvement that may be achieved by performing topological changes during the mesh optimisation process, although in practice it would be necessary to restrict the reduction in the number of vertices. Although not the subject of this paper, it is intended to add topological changes to Mesquite at a later date.
\vspace{20mm}
\begin{table}
\centering
\caption{Highest Quality Achievable}
\label{tab:1}       
\begin{tabular}{llc}
\hline\noalign{\smallskip}
Mesh & Smoother & Min/Max Dihedral Angle\\
\noalign{\smallskip}\hline\noalign{\smallskip}
\textsc{Dragon} & Inverse-Mean & 34.3$^\circ$-114.5$^\circ$ \\
& Log-Barrier & 40.9$^\circ$-107.2$^\circ$ \\
& Stellar & 40.3$^\circ$-118.9$^\circ$ \\
\textsc{Concrete} & Inverse-Mean & 26.1$^\circ$ 138.77$^\circ$ \\
 & Log-Barrier &30.94$^\circ$-129.82$^\circ$ \\
 & Stellar & 28.49$^\circ$-139.97$^\circ$ \\
\textsc{Graphite} & Inverse-Mean & 20.536$^\circ$-148.12$^\circ$ \\
& Log-Barrier & 26.103$^\circ$-138.77$^\circ$ \\
& Stellar & 22.4$^\circ$-148.83$^\circ$ \\
\noalign{\smallskip}\hline
\end{tabular}
  \label{tab:1} 
\end{table}


  %
  %

\subsection{Patch Improvement}
\label{subsubsec:2.1.1}
As we are focussed on improving the worst elements in a mesh, it is inefficient to operate on all elements. Therefore, a modified form of patch-based improvement is used. Patch improvement involves breaking the mesh up into smaller mesh patches and improving each patch individually. Since we only wish to operate on the worst elements, we select the patches containing these elements and improve them. A target quality is selected and all elements with quality worse than this target are used to generate the patches. This is similar to the method employed by \cite{Alliez2005}. It has been found that the most efficient approach is to adopt an iterative process, looping over successive patches. Results presented in Section \ref{results} show just how effective this method is at reducing the time taken to improve meshes.

\section{Surface Mesh Optimisation}
\label{sec:3}
If the worst element of a mesh lies on a boundary, then it becomes difficult to improve the mesh without changing the geometry of the domain. Several methods have been developed to tackle this \cite{Knupp2012}. If the domain of the boundary is a straight line or a planar surface there are two possible options. Mesquite provides built in functionality which "snaps" nodes which have been moved from either a planar surface or a straight line back onto the correct domain. This method is effective for surfaces which may be mathematically defined, but is not sufficient for more complex ones. 
\subsection{Surface Quadrics}
Klingner~\cite{Klingner2008} developed a method of using surface quadrics and implemented this into Stellar. This method assigns an error to a vertex which has been moved based on how far it has moved from the planes created by the original triangular faces that adjoined it \cite{Klingner2008}. This approach is summarised here. Let \emph{P} be the set of planes created by the surface triangular faces adjoining a vertex, \emph{v}. The quadric error for a point \emph{x} relative to \emph{v} is defined as 

\begin{equation}
  \label {eqn:quad}
  Q_{v}(x)=\Sigma \delta_{i}(x)^{2}
\end{equation}
where \textit{$\delta$}$_{i}(x)$ is the perpendicular distance of \emph{x} from the \textit{i}$^{th}$ plane. This means that if a vertex moves along a surface, there is no quadric error. However, if a vertex moves perpendicular to a surface, the quadric error increases rapidly. By limiting the quadric error, the amount by which a vertex may move from a surface is limited~\cite{Klingner2008}. A penalty function is used to trade the quality of an element off against its quadric error. Klingner~\cite{Klingner2008} has shown that it is possible to achieve high quality improvement by making small changes to the surface of a mesh.

Although using surface quadrics has been shown to be effective, this method has the disadvantage that the geometry of the domain is being changed and there is no guarantee that the volume will remain constant. Therefore, we wish to develop a method whereby surface vertex movement does not change the geometry of the domain, using only information which may be derived from the discretised domain.

\subsection{Optimising Mesh Surface Using Boundary Representation }
All modern CAD systems use a Boundary Representation or B-Rep solid model to store geometry. If this information is available, then Mesquite can optimise surface meshes whilst restricting surface vertices to their respective surfaces. Mesquite also contains a deforming domain class whereby the initial mesh of the undeformed domain is used to guide the optimisation of deformed mesh, \cite{Knupp2012}. However, there are many cases of domain deformation whereby this information will not be available, for example crack propagation and problems involving free fluid surfaces, e.g. dam break or microfluids with surface tension where the evolution of geometry is unknown and governed by physical equations. No additional information  about such surfaces is available so any optimisation of the surface must be based on information extracted from the discretised geometry. 

\subsection{Generating Surface Constraints from the Discretised Domain}\label{surface}

This section discusses the development and implementation of an algorithm which allows for the movement of nodes on a non-planar surface. This algorithm does not change the underlying mesh geometry as it is based on our hypothesis that for a given shape, the volume to surface area ratio is a constant. 
\begin{equation}
  \frac{V}{A}=C_{0}
  \label{six}
\end{equation}
\noindent where \emph{V} is the domain volume, \emph{A} is the surface area of the domain and \emph{C} is a constant.\\

\noindent From (\ref{six})
\begin{equation}
\label{folle}
\int_V \textrm{d}V = C_0 \int_A \textrm{d}A 
\end{equation}
Using the divergence theorem, the volume integral becomes a surface integral:
\begin{equation}
\label{fou}
\int_V \textrm{d}V = 
\frac{1}{3} 
\int_V  \textrm{div}(\mathbf{X}) \textrm{d}V =
\frac{1}{3} 
\int_A \mathbf{X}\cdot \frac{1}{\| \mathbf{N} \|}\mathbf{N} \textrm{d}A 
\end{equation}
Where  $\mathbf{X}$ is a Cartesian coordinate and $\mathbf{N}$ is the outward pointing normal of the surface. Combining equations (\ref{folle}) and (\ref{fou}):
\begin{equation}
\frac{1}{3} 
\int_A \mathbf{X}\cdot\frac{1}{\| \mathbf{N} \|}\mathbf{N} \textrm{d}A 
=
C_0 \int_A \textrm{d}A 
\end{equation}
Rewriting the above gives:
\begin{equation}
\nonumber
\frac{1}{3} 
\int_A (\mathbf{X}\cdot\frac{1}{\| \mathbf{N} \|}\mathbf{N}-C_1) \textrm{d}A=0 \;\;  \mathrm{where}\;\;C_1 = 3\,C_0.
\end{equation}
which yields a local variant as follows:
\begin{equation}
\mathbf{X}\cdot\frac{\mathbf{N}}{\| \mathbf{N} \|} = C_1
\end{equation}
A first order Taylor Series yields:
\begin{equation}
\nonumber
\mathbf{X}_i\cdot\frac{\mathbf{N}_i}{\| \mathbf{N}_i \|} +
\frac{\mathbf{N_i}}{\| \mathbf{N}_i \|}\cdot\frac{\partial \mathbf{X}_i}{\partial \mathbf{X}_i}\delta\mathbf{X}_{i+1} +
\mathbf{X_i}\cdot \frac{1}{\| \mathbf{N}_i \|}\frac{\partial \mathbf{N}_i}{\partial \mathbf{X}_i}\delta\mathbf{X}_{i+1} -\\
\end{equation}
\begin{equation}
(\mathbf{X}_i\cdot\mathbf{N}_i) \frac{\mathbf{N}_i}{\| \mathbf{N}_i \|^{3}} 
\frac{\partial \mathbf{N}_i}{\partial \mathbf{X}_i}\delta\mathbf{X}_{i+1} 
= C_1 
\end{equation}
Where $\delta$ represents an iterative change. Rearranging,
\begin{equation}
\nonumber
\frac{\mathbf{N_i}}{\| \mathbf{N}_i \|} \cdot \delta\mathbf{X}_{i+1} +
\mathbf{X_i} \cdot \frac{1}{\| \mathbf{N}_i \|}\frac{\partial \mathbf{N}_i}{\partial \mathbf{X}_i}\delta\mathbf{X}_{i+1} -
\end{equation}
\begin{equation}
(\mathbf{X}_i\cdot\mathbf{N}_i) \frac{\mathbf{N}_i}{\| \mathbf{N}_i \|^{3}} 
\frac{\partial \mathbf{N}_i}{\partial \mathbf{X}_i}\delta\mathbf{X}_{i+1} 
= C_1 - \mathbf{X}_i\cdot\frac{\mathbf{N}_i}{\| \mathbf{N}_i \|}
\label{integral}
\end{equation}

\noindent The second and third terms of the left hand side cancel out, leading to the following surface constraint equation:

\begin{equation}
\frac{\mathbf{N_i}}{\| \mathbf{N}_i \|} \cdot \delta\mathbf{X}_{i+1} 
= C_1 - \mathbf{X}_i\cdot\frac{\mathbf{N}_i}{\| \mathbf{N}_i \|}
\label{integral}
\end{equation}

\noindent Enforcing this equation in a weighted residual sense leads to:
\begin{equation}
  \mathbf{C}\delta\mathbf{X} = \mathbf{g}
  \label{fif}
\end{equation}

\noindent where $\mathbf{C}$ is a constraint matrix and $\mathbf{g}$ is a residual vector. This constraint equation ensures that the volume to surface area ratio is conserved.


Following \cite{Ainsworth2001}, the non-linear system of equations in (\ref{Sdx}) are modified to explictly account for the constraint equation (\ref{fif}) as follows:
\begin{equation}
  \mathbf{S'} \delta X = - \mathbf{f'}
\end{equation}
where
\begin{equation}
  \mathbf{S}'=\mathbf{C}^{T}\mathbf{C}+\mathbf{Q}^{T}\mathbf{S}\mathbf{Q}
\end{equation}
\begin{equation}
  \mathbf{f}'=\mathbf{C}^{T}\mathbf{g}+\mathbf{Q}^{T}(\mathbf{f}-\mathbf{SRg})
\end{equation}
\begin{equation}
  \mathbf{R}=\mathbf{C}^{T}(\mathbf{CC}^{T})^{-1}
\end{equation}
\begin{equation}
  \mathbf{Q}=\mathbf{I}-\mathbf{C}^{T}(\mathbf{CC}^{T})^{-1}\mathbf{C}
\end{equation}

\begin{figure*}
  \centering
 \hspace{-20mm}\begin{subfigure}[b]{0.3\textwidth}
    \center{\includegraphics[height=6cm]{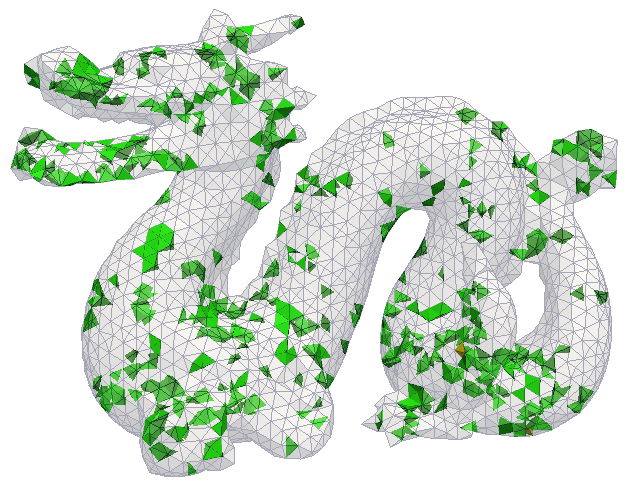}}
    \raggedright\caption{\textsc{Dragon}, 32959 Tetrahedra \cite{Klingner:2007:ATM}}
    \label{dragon_best}
  \end{subfigure}
  \hspace{35mm}
 \begin{subfigure}[b]{0.3\textwidth}
    \center{\includegraphics[height=6cm]{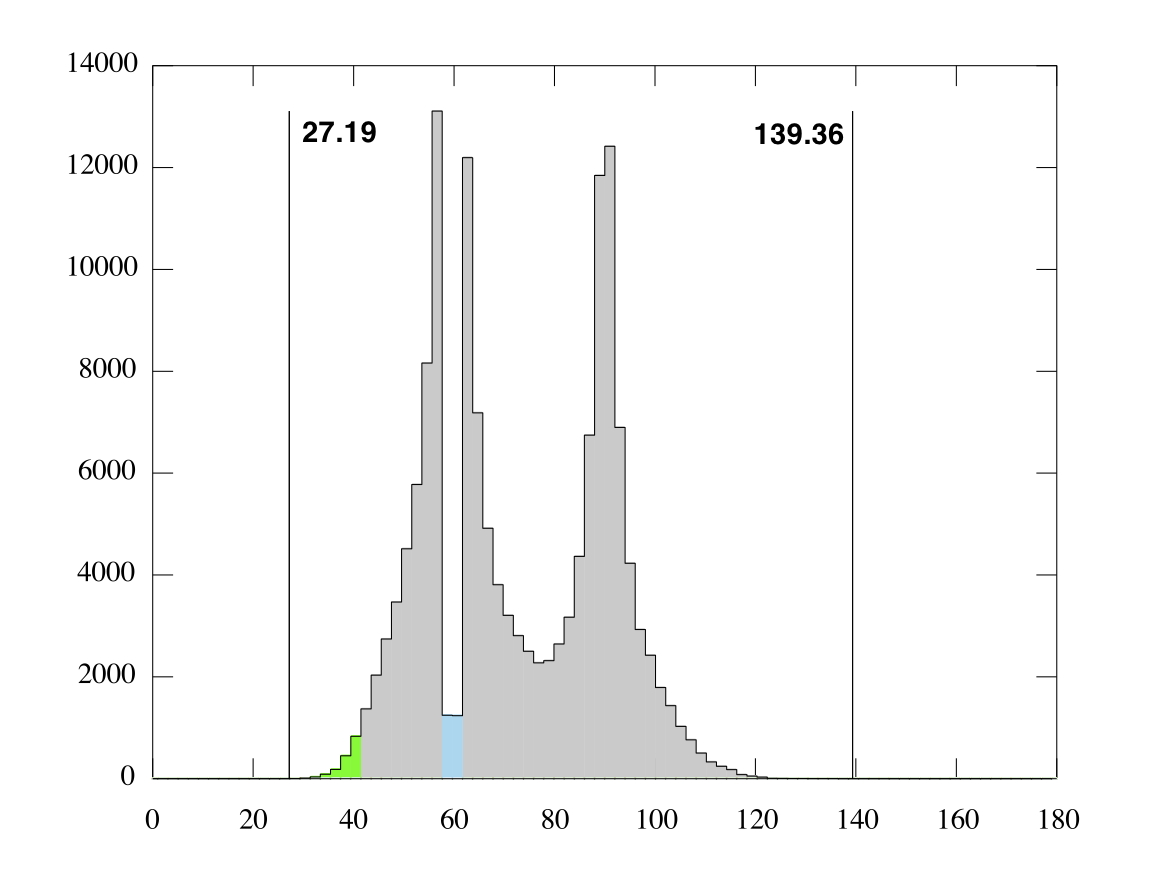}}
    \label{dragon_best_hist}
  \end{subfigure}\\

  \hspace{-20mm}\begin{subfigure}[b]{0.35\textwidth}
    \center{\includegraphics[height=6cm]{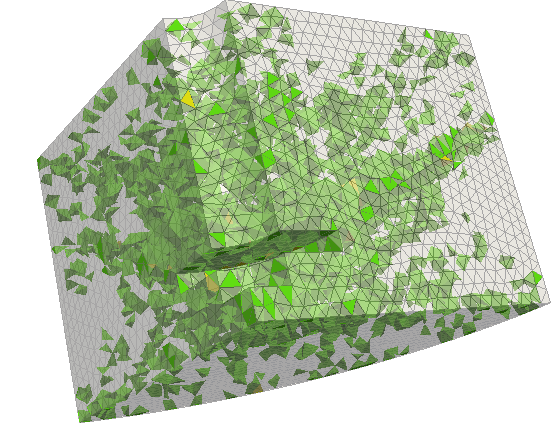}}
    \raggedleft\caption{\textsc{Concrete Cylinder}, 73864 Tetrahedra}
    \label{cow_best}
  \end{subfigure}
  \hspace{35mm}
 \begin{subfigure}[b]{0.3\textwidth}
    \center{\includegraphics[height=6cm]{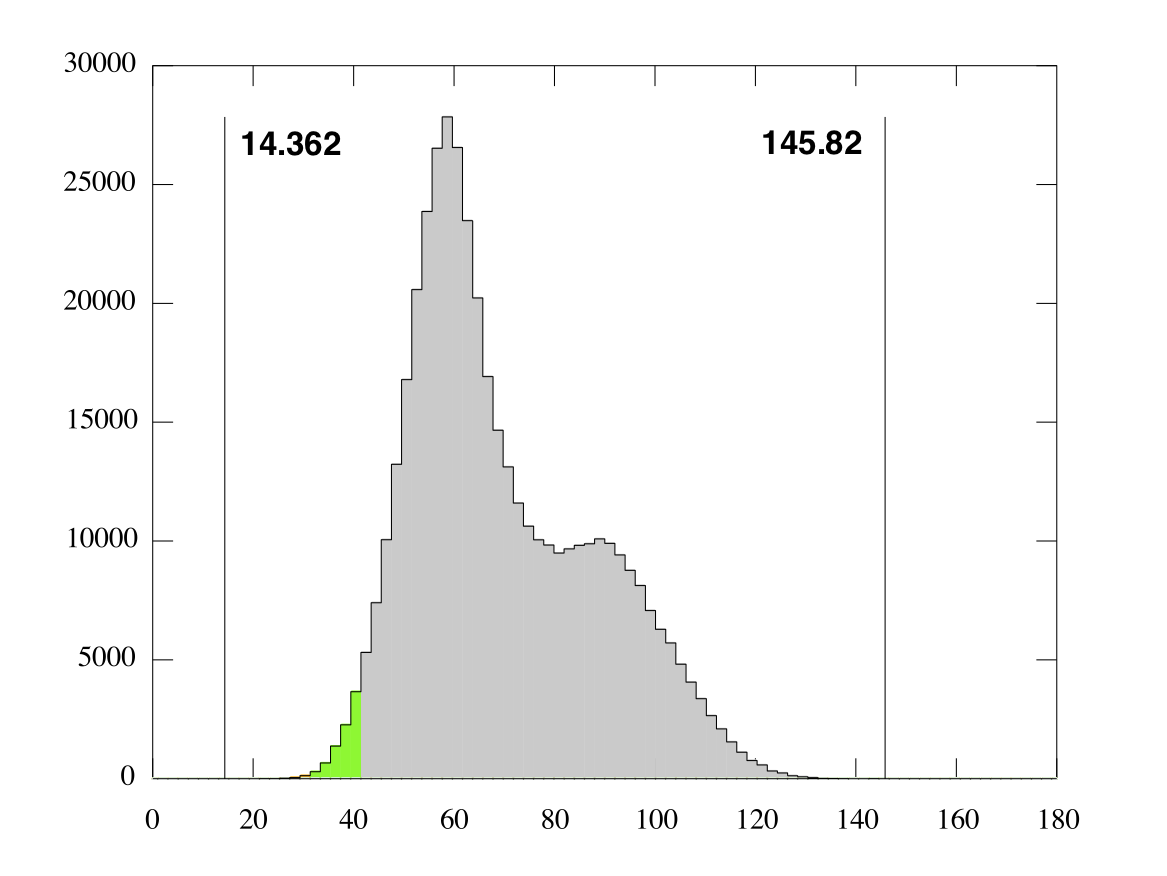}}
    \label{best_cow_hist}
  \end{subfigure}\\
  \hspace{-20mm}\begin{subfigure}[b]{0.35\textwidth}
    \center{\includegraphics[height=6cm]{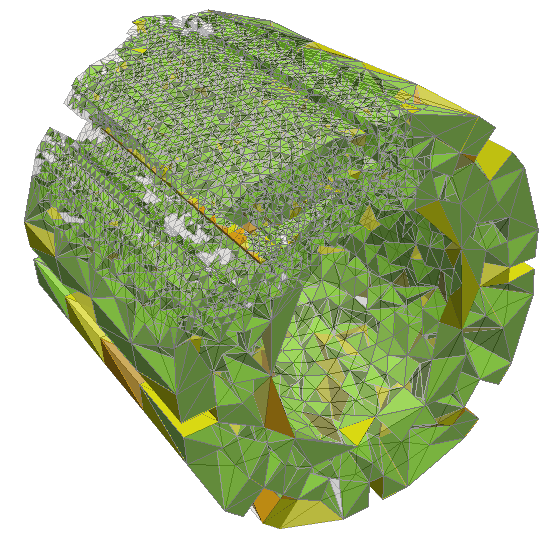}}
    \raggedright\caption{\textsc{Graphite Brick}, 100556 Tetrahedra}
    \label{stay_puft_best}
  \end{subfigure}
  \hspace{35mm}  
  \begin{subfigure}[b]{0.3\textwidth}
    \center{\includegraphics[height=6cm]{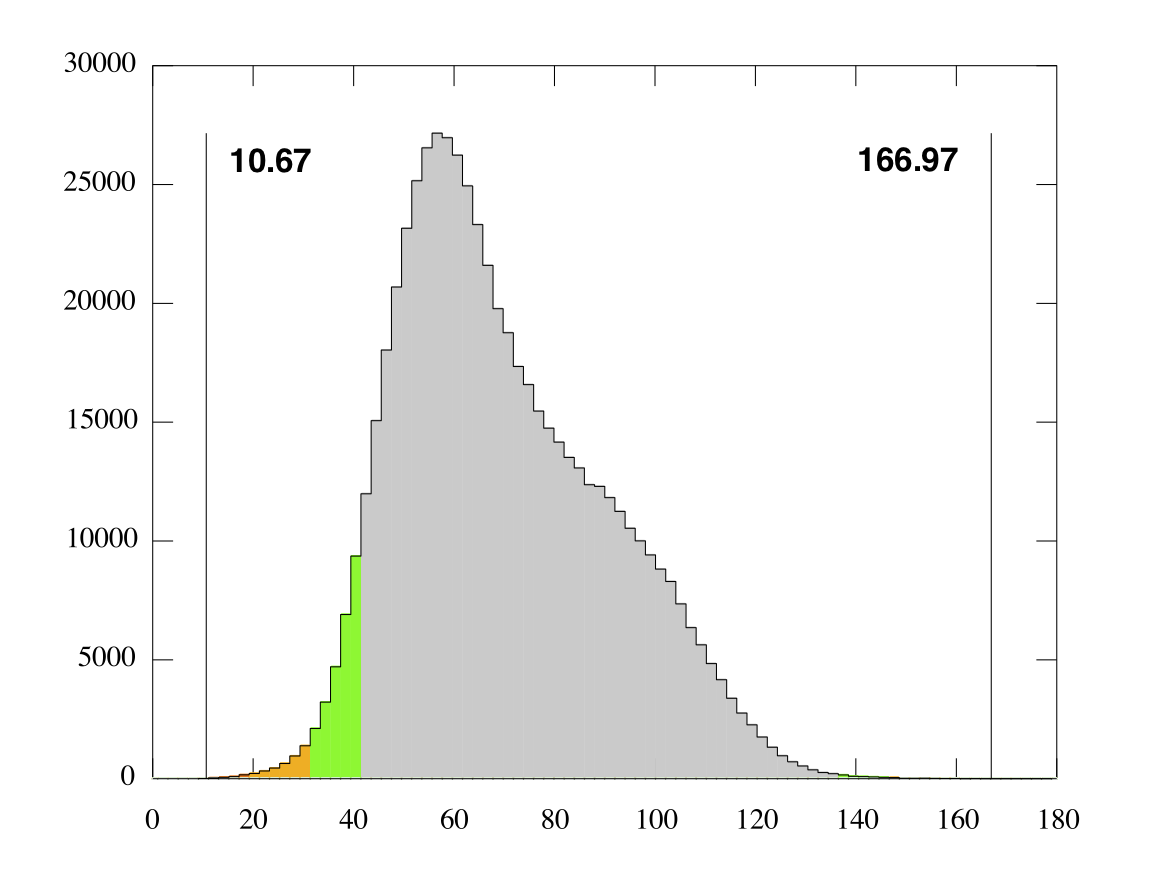}}
    \label{brick_best_hist}
  \end{subfigure}

  \caption{Results from combination of Log-Barrier, Patch-Based Improvement and Surface Optimisation techniques. Mesh used for Testing $($the height of blue columns in the histogram have been divided by 20 due to the many occurrences of these angles$)$. Orange tetrahedra have angles under 20$^\circ$ or greater than 160$^\circ$, yellow tetrahedra have angles between 20$^\circ$ and 30$^\circ$ or 150$^\circ$ and 160$^\circ$ and green tetrahedra have angles between 30$^\circ$ and 40$^\circ$ or 140$^\circ$ and 150$^\circ$.}
  \label{plot:results}
  %
  %
\end{figure*}

\section{Results and Discussion}
\label{results}
The three meshes presented in Figure \ref{plot:mesh} were optimised using the techniques described in the previous sections. The results are presented in Figure \ref{plot:results}. The worst angles in all three meshes have been eliminated leaving meshes which would be much more suitable for FE simulations. Further investigations shows that the worst elements are ones where all four nodes lie on the surface in locations of high curvature. This means that these nodes are defining the surface and any movement of them will results in unacceptable changes to the mesh geometry and changes to the volume. In situations such as this, mesh smoothing is fundamentally limited as the only way to eliminate these poor dihedral angles is to change the mesh topology. This demonstrates that the quality which may be achieved by smoothing alone is limited by the initial mesh configuration. Efforts during the mesh generation stage to ensure that all the nodes of an individual element do not lie on the surface are crucial.
The histograms in Figure \ref{plot:mesh} show that selective patch improvement has the potential to be very effective as relatively very few poor angles exist in all three mesh and this is clearly visible in the timings in Table \ref{tab:patch}. As stated in Section \ref{subsubsec:2.1.1}, a major requirement of a useful mesh optimisation toolkit is efficiency. As may be seen in Table \ref{tab:patch}, selective patch improvement greatly reduces the time taken to optimise a mesh. However, this requires that only relatively few poor quality elements exist in the mesh. For meshes with many poor quality elements, the advantages of selective patch improvement are not as great.

\begin{table}
\centering
\caption{Time Taken for Optimisation}
\begin{tabular}{lcc}
\hline\noalign{\smallskip}
& Time (s)&\\
Mesh & All Patches & Selective Patches\\
\noalign{\smallskip}\hline\noalign{\smallskip}
\textsc{Dragon} & 61.9 & 13.1\\
\textsc{Concrete Cylinder} & 265.9 & 57,4\\
\textsc{Graphite Brick} & 337.9 & 71.4 \\
\noalign{\smallskip}\hline
\end{tabular}
  \label{tab:patch} 
\end{table}

\begin{figure}
    \includegraphics[height=5cm]{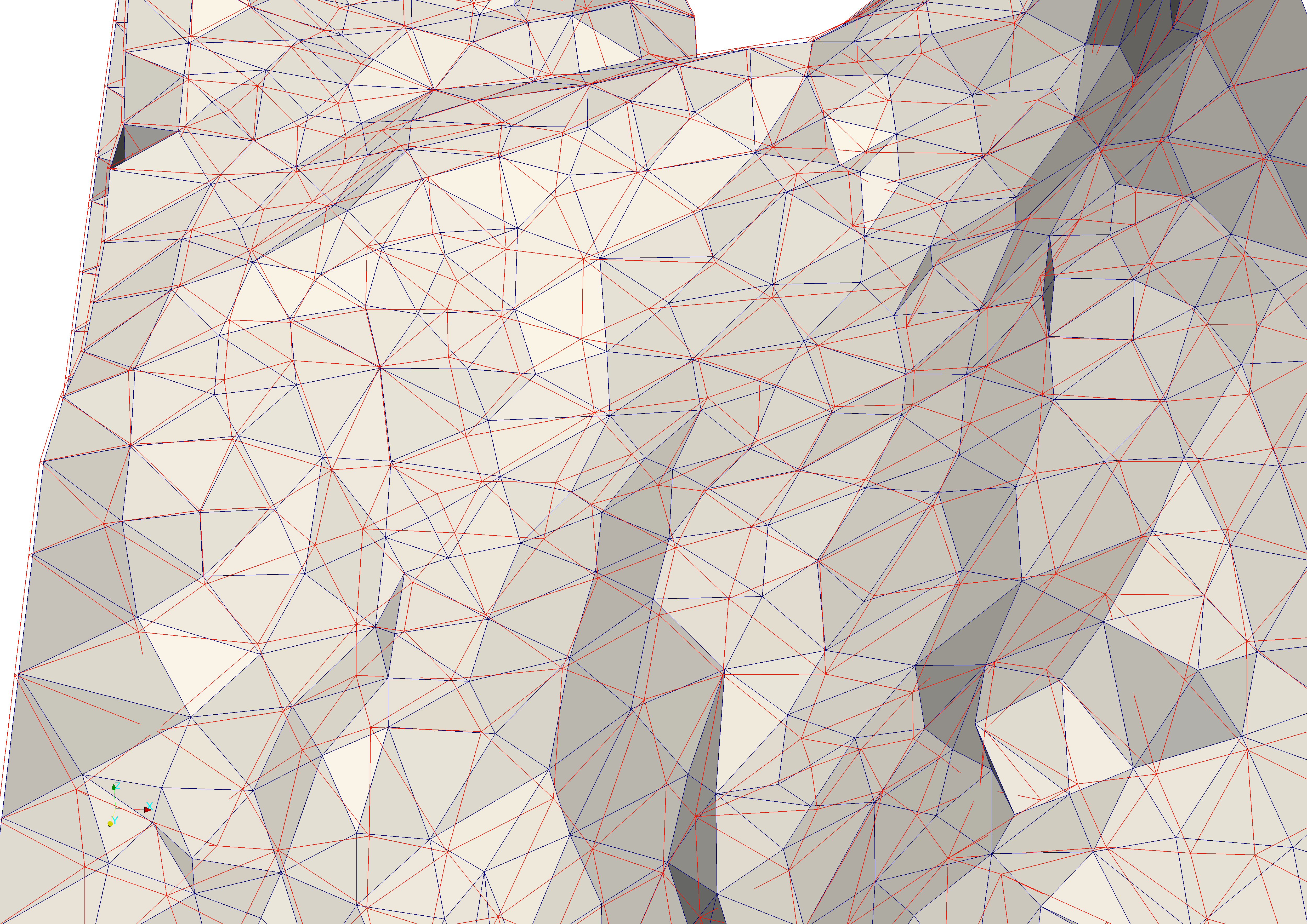}
  \caption{Effectiveness of surface mesh improvement: red mesh is optimised and blue is original mesh}
  %
  %
  \label{cow:zoom}
\end{figure}

To demonstrate the effectiveness of the surface optimisation algorithm, the original mesh and optimised mesh of \textsc{Graphite Brick} are overlain on each other so that it is possible to see the movement of surface nodes, Figure \ref{cow:zoom}. Although large nodal displacements are observed in some places, the overall crack shape is preserved. As is often the case with complex meshes, many of the worst elements have all four nodes on the surface. This means that traditional optimisation-based smoothers will not be able to improve it due to the inability to move surface nodes. The results obtained from using this complex mesh demonstrate how effective this is. Along with the preservation of mesh geometry, the volume is also completely preserved, meaning this technique may be applied to many complex simulations. For the first time, complex mesh surfaces may now be optimised whilst preserving both geometry and volume using only the mesh to define the surface. 

In this example, mesh optimisation was achieved with a change in volume of 0.0071\%, which is obviously negligible. It is also worth noting that the mesh optimisation routines included in the release version of Mesquite could not improve these meshes due to all the nodes of the worst elements being on the mesh surface. The quality improvement achieved using such complex meshes demonstrates just how effective the combination of the algorithms presented in this paper are, with the movement of surface nodes enabling the other algorithms to improve the meshes.

This method is also effective when applied to planar surfaces as can be seen in \textsc{Concrete Cylinder} and \textsc{Graphite Brick}, both of which have such surfaces.

\section{Conclusions and Future Work}
We have developed and implemented a very effective mesh optimisation methodology. The combination of a log-barrier objective function, selective-patch based improvement and surface optimisation has enabled us to optimise mesh in a robust, reliable and efficient manner which previously would not have been possible, as demonstrated by the improvement achieved in all three test meshes. The optimisation of surface nodes is completely automated and integrated into Mesquite as are the log-barrier objective function and worst patch selector. The preservation of quality during large scale deformation makes this work very useful in many different simulations. With all these features now added to Mesquite, working robustly and efficiently, it is intended to apply these to many complex FE simulations. It is also intended to make all code available to interested parties. It is clear from the tests performed using Stellar that the ability to modify the mesh topology can greatly improve a mesh and it is intended to add limited topological transformation functionality to Mesquite.


\bibliographystyle{unsrt}      
\bibliography{library}   

%
%


\end{document}